\documentclass[12pt]{article}
\usepackage{amssymb}
\title{Tournament Sequences and Meeussen Sequences}
\author{
Matthew Cook\\
Computational and Neural Systems Program\\
California Institute of Technology\\
Pasadena, CA 91125\\
{\tt cook@paradise.caltech.edu}
\and 
Michael Kleber\thanks{Partially supported by an NSF Mathematical
Sciences Postdoctoral Research Fellowship}\\
Department of Mathematics\\
Massachusetts Institute of Technology\\
Cambridge, MA 02139\\
{\tt kleber@math.mit.edu}
}
\date{Submitted: March 22, 2000; Accepted: September 5, 2000}
\newcommand{\phito}{\stackrel{\phi}{\mapsto}}  % Extra {}s mess up spacing!
\newcommand{\pf}{\noindent{\bf Proof:\,}\,}
\newcommand{\qed}{\hfill $\blacksquare$ \bigskip}
\newcommand{\tm}{{\raisebox{.75ex}{\sc\small tm}}}

\newtheorem{thm}{Theorem}
\newtheorem{defn}[thm]{Definition}
\newtheorem{prop}[thm]{Proposition}
\newtheorem{cor}[thm]{Corollary}

\newtheorem{lemma}[thm]{Lemma}
\newtheorem{recur}{Recurrence}
\begin{document}
\pagestyle{myheadings}
\markright{\rm{\sc the electronic journal of combinatorics} \textbf{7} (2000), \#R44\hfill} 
\thispagestyle{empty} 
\maketitle

\begin{abstract}
A {\em tournament sequence} is an increasing sequence of positive
integers $(t_1,t_2,\ldots)$ such that $t_1=1$ and $t_{i+1} \leq 2t_i$.
A {\em Meeussen sequence} is an increasing sequence of positive
integers $(m_1,m_2,\ldots)$ such that $m_1=1$, every nonnegative
integer is the sum of a subset of the $\{m_i\}$, and each integer
$m_i-1$ is the sum of a unique such subset.

We show that these two properties are isomorphic.  That is, we present
a bijection between tournament and Meeussen sequences which respects
the natural tree structure on each set.  We also present an efficient
technique for counting the number of tournament sequences of length
$n$, and discuss the asymptotic growth of this number.  The counting
technique we introduce is suitable for application to other
well-behaved counting problems of the same sort where a closed form or
generating function cannot be found.

\end{abstract}

\section{Introduction}
\label{sec_intro}

An infinite {\em tournament sequence} $T$ is an infinite sequence of
positive integers $T=(t_1,t_2,\ldots)$ such that
\begin{itemize}
  \item $t_1=1$ and $t_i < t_{i+1} \leq 2t_i$ for $i=1,2,\ldots$
\end{itemize}
For example, the first infinite tournament sequence in lexicographic
order is $t_i=i$, and the last is $t_i=2^{i-1}$.  A finite tournament
sequence $T=(t_1,\ldots,t_n)$ is a truncated infinite tournament
sequence.

An infinite {\em Meeussen sequence} $M$ is an infinite sequence of
positive integers $M=(m_1,m_2,\ldots)$ such that
\begin{itemize} \addtolength{\itemsep}{-6pt}
  \item $m_1=1$ and $m_i<m_{i+1}$ for $i=1,2,\ldots$,
  \item Every nonnegative integer is the sum of a subset of 
        the $\{m_i\}$, and 
  \item Each integer $m_i-1$ is the sum of a {\em unique} subset 
        of the $\{m_i\}$.
\end{itemize}
For example, the first infinite Meeussen sequence in lexicographic
order is $m_i=f_{i+1}$, the $(i+1)$st Fibonacci number, and the last
is $m_i=2^{i-1}$.  A finite Meeussen sequence $M=(m_1,\ldots,m_n)$ is
a truncated infinite Meeussen sequence.  We will see that this is
equivalent to requiring that every integer between $1$ and
$\sum_{i=1}^n m_i$ is the sum of a subset of the $\{m_i\}$.

We present a bijection $\{T\}\leftrightarrow\{M\}$ between these two
types of sequences.  The bijection is defined in
Section~\ref{sec_iso}; it preserves the length of the sequence and
respects lexicographic ordering.  It also acts in a surprising way on
sequences with certain recurrence relations, as discussed in
Section~\ref{sec_fib}.

Counting finite tournament (or equivalently Meeussen) sequences of
length $n$ is straightforward (\cite{EIS}, sequence A008934), but if
done in the obvious way takes time exponential in $n$.  In
Section~\ref{sec_count} we present an efficient polynomial-time
algorithm for producing the numbers.  The technique is suitable for
application to other well-behaved counting problems of the same sort
where a closed form or generating function cannot be found.  We also
discuss the asymptotic growth, proving that the $\log_2$ of the number 
of sequences of length $n$ is ${n \choose 2} - \log_2(n!) + O(\log(n)^2)$.

Finite tournament sequences were studied in \cite{knockout} under the
name ``random knock-out tournaments.''  A sequence $(t_1,\ldots,t_n)$
represented a tournament of $n$ rounds beginning with $1+\sum t_i$
players; in the first round $2t_n$ players are paired off randomly and
the $t_n$ losers are eliminated, leaving a tournament corresponding to
$(t_1,\ldots,t_{n-1})$.  The paper concerns the probabilities of
certain pairings occurring in such a tournament.

The count which first appeared in \cite{EIS} was performed by
M.~Torelli, who found the notion of a tournament sequence useful in
his investigation of sequences with certain properties relating to
Goldbach's conjecture~\cite{T}.  Tournament sequences also appear
independently in the work of J.~Shallit, where they are the possible
subword complexities of infinite non-periodic bit
strings~\cite{subword}.

The observation that the beheaded Fibonacci sequence satisfies the
property claimed above was made by Wouter Meeussen [private
communication, 1999], and we here name sequences with this property
Meeussen sequences in his honor.  Part of their definition is similar
to that of so-called regular sequences~\cite{regular}, in which each
term is a partial sum of preceding terms, which arise in the study of
finite probability measures.  The unique representability condition is
reminiscent of $1$-additive sequences (see \cite{additive}, for
example), but the precise form of the condition seems new.

The authors would like to acknowledge W.~Meeussen for suggesting the
question, and N.~J.~A.~Sloane for his Encyclopedia of Integer
Sequences~\cite{EIS}, which led us to notice the coincidence.  Thanks
also to J.~Polito for useful conversations, J.~Shallit for helpful
comments on an earlier draft of this work, and D.~Knuth for excellent
suggestions about the asymptotics questions.

\section{An Isomorphism on Trees}
\label{sec_iso}

\begin{figure}
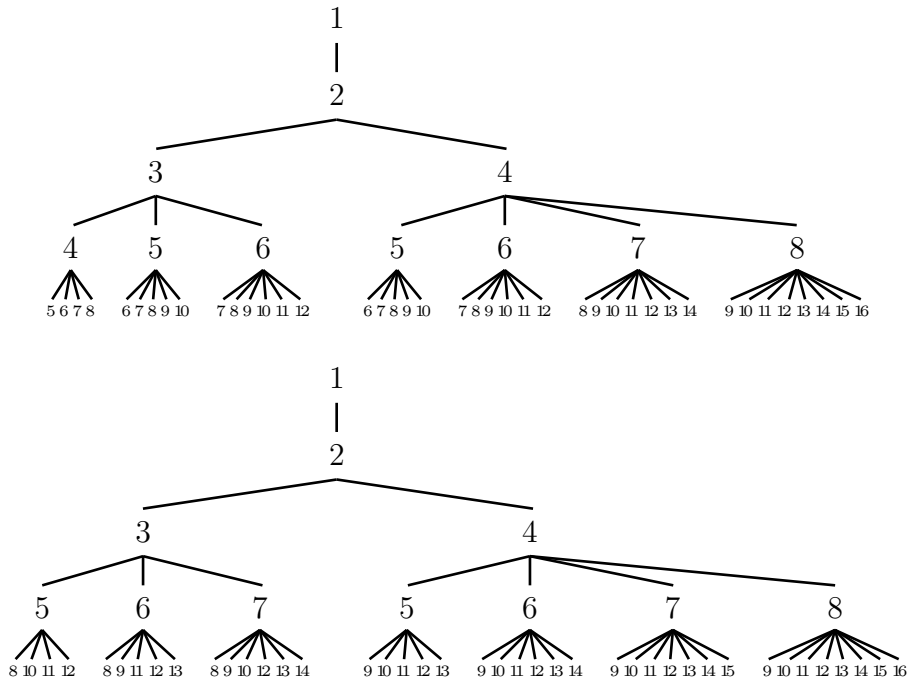

%
%  Post-script \special hacking like this is AWFUL latex style.
%  NEVER DO THIS!  --mk
%
\newcommand{\sss}{\scriptscriptstyle}

\begin{center}
$
\begin{array}{cccccccc}
 & && \hspace{-1em}1\hspace{-1em} \\
 & && \special{" 0 br}\\
 & && \hspace{-1em}2\hspace{-1em} \\
 & && \special{" mark -68 64 tr} \\
 & 3 && & & 4 \\
& \special{" mark -31 0 40 tr} &&&& \special{" mark -39 0 51 110 tr} \\
4 & 5 & 6 && 5 & 6 & 7 & 8 \\[-1pt]
\special{" mark -7 -2 3 8 tr} 
& \special{" mark -11 -6 -1 4 9 tr}
& \special{" mark -15 -10 -5 1 7 14 tr}
&& \special{" mark -11 -6 -1 4 9 tr}
& \special{" mark -16 -10 -5 1 7 14 tr}
& \special{" mark -20 -14 -8 -2 5 12 19 tr}
& \special{" mark -25 -18 -11 -3 3 10 17 24 tr}
\\[-6pt]
{\sss 5\,6\,7\,8} 
& {\sss 6\,7\,8\,9\,1\!0} 
& {\sss 7\,8\,9\,1\!0\,1\!1\,1\!2} 
&& {\sss 6\,7\,8\,9\,1\!0} 
& {\sss 7\,8\,9\,1\!0\,1\!1\,1\!2} 
& {\sss 8\,9\,1\!0\,1\!1\,1\!2\,1\!3\,1\!4} 
& {\sss 9\,1\!0\,1\!1\,1\!2\,1\!3\,1\!4\,1\!5\,1\!6}
\end{array}
$
\\[1em]
$
\begin{array}{cccccccc}
 & && \hspace{-1em}1\hspace{-1em} \\
 & && \special{" 0 br}\\
 & && \hspace{-1em}2\hspace{-1em} \\
 & && \special{" mark -73 74 tr} \\
 & 3 && & & 4 \\
& \special{" mark -38 0 45 tr} &&&& \special{" mark -46 0 54 115 tr} \\
5 & 6 & 7 && 5 & 6 & 7 & 8 \\[-1pt]
\special{" mark -10 -4 3 10 tr} 
& \special{" mark -14 -8 -2 5 12 tr}
& \special{" mark -17 -11 -5 2 9 16 tr}
&& \special{" mark -14 -8 -1 6 13 tr}
& \special{" mark -18 -11 -4 2 9 16 tr}
& \special{" mark -21 -14 -7 -1 6 13 21 tr}
& \special{" mark -25 -17 -10 -3 3 10 17 24 tr}
\\[-6pt]
{\sss 8\,1\!0\,1\!1\,1\!2} 
& {\sss 8\,9\,1\!1\,1\!2\,1\!3} 
& {\sss 8\,9\,1\!0\,1\!2\,1\!3\,1\!4} 
&& {\sss 9\,1\!0\,1\!1\,1\!2\,1\!3} 
& {\sss 9\,1\!0\,1\!1\,1\!2\,1\!3\,1\!4} 
& {\sss 9\,1\!0\,1\!1\,1\!2\,1\!3\,1\!4\,1\!5} 
& {\sss 9\,1\!0\,1\!1\,1\!2\,1\!3\,1\!4\,1\!5\,1\!6}
\end{array}
$
\end{center}
\caption{The top five rows of the trees of tournament and Meeussen sequences.}
\label{fig_trees}
\end{figure}

In this section we define a map which sends any tournament sequence,
finite or infinite, to a Meeussen sequence of the same length.  The
set of all sequences of either type can naturally be viewed as a
rooted tree: the nodes on level $n$ of the tree correspond to the
sequences of length $n$, and the parent of the sequence
$(s_1,\ldots,s_n)$ in the tree is the sequence $(s_1,\ldots,s_{n-1})$.
Our map is an isomorphism on the tree structures of the two types of
sequences.  Figure~\ref{fig_trees} shows how the beginnings of these
trees look; the node for $(s_1,\ldots,s_n)$ has the label $s_n$
written on it.

The definition of a tournament sequence $(t_1,\ldots,t_{n+1})$ says
that $t_{n+1}$ can have any value between $t_n+1$ and $2t_n$.
Therefore the tree of tournament sequences has a convenient local
description: the top node is labelled $1$, and any node labelled $k$
has $k$ children, with labels $k+1$, $k+2,\ldots,$ $2k$, respectively.
Trees with such local descriptions have been called {\em generating
trees}, and were introduced in~\cite{CGHK}.  They have been championed
by J.~West~(\cite{W1},\cite{W2}), who has used them to study
pattern-avoiding permutations, and by Barcucci {\em et.~al.}, who
applied them to the enumeration of combinatorial objects~\cite{bdpp};
see~\cite{mbm} and references therein for more on the subject.  In
West's transparent notation, the tree of tournament sequences is $$
\begin{array}{ll}
\mbox{Root:} & (1) \\
\mbox{Rule:} & (k) \to (k+1)(k+2)\ldots(2k)
\end{array}
$$

We will prove that the tree of Meeussen sequences is isomorphic to the
tree of tournament sequences by showing that it also has the property
that if a node has $k$ children, then those children have $k+1$,
$k+2,\ldots,$ $2k$ children, respectively.  Since this tree structure
clearly has no automorphisms, we conclude that there is a unique
bijection between the two trees, and therefore a unique bijection
between tournament sequences and Meeussen sequences which respects the
ideas of extending or truncating a sequence.

To better understand Meeussen sequences, we introduce some notation.

\begin{defn}
For $A=(a_1,a_2,\ldots)$ an integer sequence, finite or infinite, define:
\begin{enumerate}
  \item
    $r(A)$ to be the set of integers which are representable as
    $a_{i_1}+\ldots+a_{i_n}$ for some $a_{i_1},\ldots,a_{i_n}$ in $A$, 
    $i_1<\cdots<i_n$,
  and
  \item
    $ur(A)$ to be the set of integers so representable in exactly one way.
\end{enumerate}
\end{defn}
For example, if $A=(1,2,3)$ is our sequence, then
$r(A)=\{0,1,2,3,4,5,6\}$, and $ur(A)=\{0,1,2,4,5,6\}$, where $3$ is
omitted because it can be represented as $3$ and as $1+2$.  An
infinite increasing integer sequence $M=(m_1,m_2,\ldots)$ with $m_1=1$
is Meeussen if $r(M)={\mathbb Z}_{\geq0}$ and $m_i-1 \in ur(M)$ for
all $i$.

Given a finite Meeussen sequence $M=(m_1,\ldots,m_n)$ which we wish to
extend, we must certainly pick $m_{n+1}$ to be $u+1$ for some element
$u\in ur(M)$ with $u\geq m_n$.  We say such choices of $u$ are {\em
candidates}.  For example, with $M=(1,2,3)$, there are three
candidates $4,5,6$ for $m_4-1$, so $m_4$ must be one of $5,6,7$.

Our claim that the two trees are isomorphic then reduces to the
following:

\begin{prop}
\label{main}
Suppose $M=(m_1,\ldots,m_n)$ is a finite Meeussen sequence, and there
are $k$ candidates for $m_{n+1}-1$, which we designate
$u_1<u_2<\cdots<u_k$.  Then for each $j$, $1\leq j\leq k$, the
extended sequence with $m_{n+1}=u_j+1$ is also Meeussen, and has $k+j$
candidates for $m_{n+2}-1$.
\end{prop}

\pf
Let $S$ denote the sum $\sum_{i=1}^n m_i$ of the sequence; note that
$S$ is also the largest candidate, $u_k$.  Let
$M'=(m_1,\ldots,m_{n+1})$ be the extended sequence we get by choosing 
$m_{n+1}=u_j+1$, and $S'=S+m_{n+1}$ be its sum.

First, we can see by induction that $r(M')$ is the entire interval of
integers $[0,S']$.  Each representable sum in $r(M')$ is an element of
$r(M)$ or is $m_{n+1}$ added to such an element.  Assume inductively
that each of these types forms a single interval, $[0,S]$ and
$[m_{n+1},S']$, respectively.  There is no gap between the two
intervals, since $m_{n+1}-1=u_j\in r(M)$, and the induction holds.

Note that these two intervals have some overlap $[m_{n+1},S]$,
possibly empty if we chose $m_{n+1}=u_k+1=S+1$.  Anything in the
overlap can be represented in at least two ways, one with $m_{n+1}$
and one without.  Thus the candidates for $M'$, the elements of
$ur(M')$ larger than $m_{n+1}$, are in fact all larger than $S$.

Next we observe that both $r(M')$ and $ur(M')$ are invariant under the
involution $t\mapsto S'-t$, corresponding to taking the complement of
a subset.  The overlap $[m_{n+1},S]$ is similarly invariant and right
in the middle, and for each candidate $v$ of $M'$, there will be a
corresponding member $S'-v$ of $ur(M')$ smaller than $m_{n+1}$.

Similarly, $ur(M)$ has the same property: aside from $u_1,\ldots,u_k$
it contains exactly $k$ more uniquely representable numbers, all of
the form $S-u_i$.  We now know that the elements of $ur(M)$ are, in
order, $S-u_k,\ldots,S-u_1,u_1,\ldots,u_k$.

We have also concluded that the candidates for $M'$ are exactly the
numbers of the form $u+m_{n+1}$ such that $u\in ur(M)$ and the total
is strictly larger than $S$.  Therefore when we chose $m_{n+1}=u_j+1$,
the candidates for $M'$ are exactly the $k$ numbers $u_i+u_j+1$ for
$i=1,\ldots,k$ and the $j$ numbers $(S-u_i)+u_j+1$ for $i=1,\ldots,j$.
Thus there are $k+j$ candidates, as desired.
\qed

Note that we incidentally showed that half of our definition of
Meeussen sequences is unnecessary.  If we always choose the term
$m_{n+1}$ to be one more than a representable sum from $r(M)$, we
argued above that at each finite stage, $r(M)$ is the entire interval
$[0,S]$.  Thus an infinite sequence $M$ generated this way will
automatically have $r(M) = {\mathbb Z}_{\geq0}$, a condition imposed
in the original definition.

\begin{cor}
There is a unique isomorphism $\phi : \{T\} \to \{M\}$ from the set of 
tournament sequences to the set of Meeussen sequences which preserves
the rooted tree structure.  Moreover, $\phi$ respects the
lexicographic ordering on sequences of each type.
\end{cor}

\pf
By Proposition~\ref{main}, the tree structure of Meeussen sequences is
precisely that of tournament sequences: both trees start with a node
with one child, and if a node has $k$ children, then those children
have $k+1$, $k+2,\ldots,$ $2k$ children, respectively.  Since the
children of a node are all distinguishable from one another, by
virtue of their distinct numbers of children, there is a unique
bijection between the two trees.

The nodes in the trees are indexed by finite sequences of each type,
so $\phi$ is immediately defined for finite sequences.  In both trees, 
the distance from the root determines the length of the sequence, so
$\phi$ preserves it.  Infinite sequences can be thought of as infinite 
paths in the tree heading away from the root, and therefore the tree
bijection lets us define $\phi$ in this case as well.

It is evident that $\phi$ respects the lexicographic ordering since,
in the proof of Proposition~\ref{main}, we showed that when extending
a sequence $M$, the larger candidates correspond to the nodes with
more children, just as in tournament sequences.
\qed

The proof of Proposition~\ref{main} also gives us a way to calculate
$\phi(T)$ for any $T=(t_1,\ldots,t_n)$ without constructing the full
set $ur(M)$, a task that could require exponential time.  We make
repeated use of the fact that the candidates for extending
$M=(m_1,\ldots,m_n)$ are easily expressed in terms of the candidates
for $(m_1,\ldots,m_{n-1})$.  Letting $u(n,k)$ be the $k$th smallest
candidate for extending $(m_1,\ldots,m_n)$, we get the following
recurrence:
\begin{eqnarray*}
m_n &=& u(n-1,t_n-t_{n-1}) + 1,\mbox{~and} \\[4pt]
u(n,k) &=&
  \left\{ \begin{array}{ll}
    u(n-1,k-(t_n-t_{n-1})), & \mbox{if } k > t_n-t_{n-1} \\[4pt]
    S_n - u(n-1, t_n-t_{n-1}+1-k), & \mbox{if } k\leq t_n-t_{n-1}
  \end{array} \right. \\
&& \mbox{\qquad where $S_n=m_1+\cdots+m_n$.}
\end{eqnarray*}
Beginning with $m_1=1$ and $u(1,1)=1$, we can quickly calculate
each successive term of $\phi(T)$ in linear time.

\section{Properties of the Bijection}
\label{sec_fib}

In the introduction, we stated without proof that the beheaded
Fibonacci sequence $(1,2,3,5,8,13,\ldots)$ is a Meeussen sequence, and
moreover that it is the smallest one in lexicographic order, the image
of $(1,2,3,4,5,6,\ldots)$ under the map $\phi$ defined above.  This is
a special case of a more general surprising property of $\phi$, which
we prove in this section.

Since $(1,2,3,5,8,13,\ldots)$ is, coincidentally, again a tournament
sequence, we can apply $\phi$ to it as well.  This leads to the
following computational observation:
$$
(1, 2, 3, 5, 8, 13, 21,\ldots) \phito
(1, 2, 3, 6, 11, 20, 37,\ldots) \phito
(1, 2, 3, 7, 13, 25, 48,\ldots)
$$
The middle sequence is the ``3-bonacci'' sequence beginning $(1,2,3)$
and the last is the ``4-bonacci'' sequence beginning $(1,2,3,7)$,
where we say a sequence is $k$-bonacci if each term (after the first
$k$) is the sum of the previous $k$ terms.
Further experimentation reveals that, for example,
$$
(1, 2, 4, 7, 12, 20, 33, 54, 88, 143,\ldots) \phito
(1, 2, 4, 7, 13, 24, 44, 81, 149, 274,\ldots).
$$
The first sequence begins $(1,2)$ and thereafter each term is one plus
the sum of the previous two, while the second sequence is the
3-bonacci sequence beginning $(1,2,4)$, with no additive constant in
the recurrence.

Inspired by examples of this type, we make the following observation.

\begin{prop}
\label{fib}
Suppose $(t_1,\ldots,t_{n+1})$ is a finite tournament sequence, and
suppose that for some integers $k,c$, it happens that both
$$
\begin{array}{llll}
t_{n+1}  &=&  t_{n  }+t_{n-1}+\cdots+t_{n-k+1}+c & \mbox{and} \\
t_{n  }  &=&  t_{n-1}+t_{n-2}+\cdots+t_{n-k  }+c.
\end{array}
$$
Then $(t_1,\ldots,t_{n+1})\phito(m_1,\ldots,m_{n+1})$, where
$$
\begin{array}{rcl}
m_{n+1}     &=&   m_{n  }+m_{n-1}+\cdots+m_{n-k+1}+m_{n-k  }.
\end{array}
$$
\end{prop}

While the statement of the proposition is designed to mimic the
examples above, the hypothesis simplifies to $t_{n+1} =
2t_{n}-t_{n-k}$.  Since our map $\phi$ respects extending or
truncating a sequence, we are really just assuming that $t_{n+1} =
2t_{n}-t_{n-k}$ for some single pair $n,k$; the effect is purely local.

To prove the proposition, it would help to have a more concrete
relationship between terms of the tournament and Meeussen sequences
associated to one another by our bijection.

\begin{lemma}
\label{lemma_tm}
Let $T=(t_1,t_2,\ldots)$ be a (finite or infinite) tournament sequence
with associated Meeussen sequence $\phi(T)=M=(m_1,m_2,\ldots)$.  Write
$ur(M)$, the uniquely representable sums of $M$, as
$\{u_1<u_2<u_3<\cdots\}$.  Then:
\begin{enumerate}
\item
The $t_i$'th uniquely representable sum $u_{t_i}$ is $m_i-1$, and
\item
The next uniquely representable sum $u_{t_i+1}$ is 
$m_1+m_2+\cdots+m_{i-1} + 1$.
\end{enumerate}
\end{lemma}

\pf
\begin{enumerate}
\item
When $i=1$, we check that $t_1=1$, $u_1=0$, and $m_1=1$ directly.  Now
assume by induction that the map $k\mapsto1+u_k$ sends $t_i$ to $m_i$
for some $i$.  Then it sends $t_i+1,t_i+2,\ldots$ to the first,
second$,\ldots$ number in $ur(M)$ greater than $m_i$, which we showed
was the desired value during the proof of Proposition~\ref{main}.
\item
Consider the truncated sequence $M_i=(m_1,\ldots,m_i)$.  The
map $t \mapsto (m_1+\cdots+m_i)-t$, as we saw in the proof of
Proposition~\ref{main}, is an involution on $ur(M_i)$, and takes
$m_i-1$ to $m_1+\cdots+m_{i-1} + 1$.  Certainly nothing in between is
uniquely representable; this is the ``overlap'' interval of numbers
which can be represented either with or without using $m_i$.  This
in turn means that $m_{i+1}$ is strictly larger than 
$m_1+\cdots+m_{i-1} + 1$, which is therefore in $ur(M)$ since it is in 
$ur(M_i)$.
\qed
\end{enumerate}

\noindent{\bf Proof\,} of Proposition~\ref{fib}:
Consider $M=(m_1,\ldots,m_n)$ with sum $S=m_1+\cdots+m_n$ and
$ur(M)=\{u_1<u_2<u_3<\cdots\}$.  By the first part of
Lemma~\ref{lemma_tm}, we know that $m_n=1+u_{t_n}$.  Therefore $ur(M)$
contains exactly $2t_n$ numbers: $u_1<\cdots<u_{t_n}$, which are less
than $m_n$, and another $t_n$ which are their images under the
involution $t\mapsto S-t$.  In particular, $u_{2t_n-i}=S-u_{i+1}$.

Now consider what happens when we pick $t_{n+1}=2t_n-t_{n-k}$, as
supposed by Proposition~\ref{fib}, and extend $M$ accordingly:
\begin{eqnarray*}
m_{n+1} 
 &=& 1 + u_{t_{n+1}} 
                 \mbox{~by Lemma~\ref{lemma_tm},} \\
 &=& 1 + u_{2t_n-t_{n-k}} \\
 &=& 1 + S - u_{t_{n-k}+1} \\
 &=& 1 + S - (m_1+\cdots+m_{n-k-1}+1)
                 \mbox{~by Lemma~\ref{lemma_tm} again,} \\
 &=& m_n + m_{n-1} + \cdots + m_{n-k}
\end{eqnarray*}
The new term of $M$ is the sum of the previous $k$ terms, as claimed. \qed

\section{The Growth of the Tree}
\label{sec_count}

We would like to know the number of tournament or Meeussen sequences
of length $n$, which we will designate $s(n)$.  Equivalently, we want
to know the number of nodes on the $n$th level of the tree shown in
Figure~\ref{fig_trees} (p.~\pageref{fig_trees}), where we can count
that $s(n)=1,1,2,7,41$ for $n$ up to 5.  Counting the nodes directly
takes time exponential in $n$, and while we cannot present a solution
in closed form, we can offer an efficient polynomial-time algorithm.
We would also like to know the asymptotic behavior of $s(n)$ as $n$
gets large.  The asymptotics reveal that $s(n)$ grows so quickly that
its generating function cannot be algebraic.

\subsection{Exact Counting}

Throughout this section, we consider the nodes to be labelled as in the
tree of tournament sequences: a node with label $(k)$ has $k$
children, with labels $(k+1)$, $(k+2),\!$ \ldots, $(2k)$,
respectively.  Based on this definition, we note that the function
$c(n,k)$ counting the number of nodes with label $(k)$ in row $n$ of
the tree satisfies:
\begin{recur}
\label{rec_c}
\begin{eqnarray*}
c(1,k) &=& \delta_{k,1} \\
c(n,k) &=& \sum_{j=\lceil \frac{k}{2} \rceil}^{k-1} c(n-1,j) \,\,
  \mbox{\quad for $n>1$.}
\end{eqnarray*}
\end{recur}
We could then find the number of nodes on row $n$ by summing $c(n,k)$
for all $k$ up to $2^{n-1}$.  The work involved grows exponentially in
$n$, though, so for large $n$ this is impractical.  We could define a
generating function in two independent variables
$$
g(x,y) = \sum_{n,k} x^n y^k c(n,k) = xy + x^2y^2 + x^3(y^3+y^4) + \cdots
$$
which, based on Recurrence~\ref{rec_c}, must satisfy
$$
g(x,y) = xy + \frac{xy}{1-y}\, g(x,y) - \frac{xy}{1-y}\, g(x,y^2).
$$
Rewriting this as $g(x,y) = \frac{xy(y-1)}{xy+y-1} + \frac{xy}{xy+y-1}\, 
g(x,y^2)$, we can solve formally by iterated substitution to get
$$
g(x,y) = \sum_{n=0}^{\infty} \left( y^{2^n}-1 \right)
          \prod_{k=0}^n \frac{x y^{2^k}}{x y^{2^k} + y^{2^k} - 1}.
$$
This seems to offer dim prospects for a nice form for $s(n)$, the
coefficient of $x^n$ at $y=1$.

Alternatively, one could hope to work with the function $d(n,k)$ which
counts the number of $n$th-generation descendents of a node labelled
$(k)$.  We can count these descendents by summing the number of
$(n-1)$-generation descendents of the node's $k$ children:
\begin{recur}
\label{rec_d}
\begin{eqnarray*}
d(1,k) &=& k \mbox{\quad for all $k\geq1$,} \\
d(n,0) &=& 0 \mbox{\quad for all $n\geq1$, and otherwise,} \\
d(n,k) &=& \sum_{j=k+1}^{2k} d(n-1,j).
\end{eqnarray*}
\end{recur}
The number of nodes on row $n$ of our tree is then $s(n)=d(n-1,1)$.
Torelli points out that this recurrence can be expressed in closed
form, by replacing the last line with
$$
d(n,k) \,=\, d(n,k-1) - d(n-1,k) + d(n-1,2k-1) + d(n-1,2k)
$$
This alternate version embodies the notion that the tree below a node
labelled $(k)$ looks just like the tree below a $(k-1)$, but modified
by pruning the branch beginning with the child $(k)$ and grafting on
branches beginning with $(2k-1)$ and $(2k)$ instead.  However, as $n$
increases, either version still involves calculating an exponentially
growing set of values.

We offer instead the following technique for calculating the growth of
the tree in polynomial time.  Consider the family of functions
$p_n$ for $n=1,2,3,\ldots$ such that $p_n(k)$ is the number of
$n$th-generation descendents of a node labelled $(k)$, what we called
$d(n,k)$ above.  Then in the spirit of Recurrence~\ref{rec_d}, we can
get a recurrence relation for the functions $p_n$ themselves:
\begin{recur}
\label{rec_p}
\begin{eqnarray*}
p_1(k) &=& k, \\
p_n(k) &=& \sum_{j=k+1}^{2k} p_{n-1}(j).
\end{eqnarray*}
\end{recur}
Purists would start the recurrence with $p_0(k)=1$ instead.

The key observation is that each $p_n$ is in fact a degree $n$
polynomial in $k$, which we obtain by symbolic summation of a range of
values of $p_{n-1}$.  Recall that the sum $\sum_{j=1}^k j^n$ is a
polynomial in $k$ of degree $n+1$.  Our recurrence states that
$p_n(k)$ is the sum of the first $2k$ values of $p_{n-1}$ minus the
sum of the first $k$ values.  Since $p_{n-1}$ is a polynomial in $k$
by induction, so is $p_n$.  The next few polynomials after $p_1=k$ are
$$
p_2 = \frac{3k^2+k}{2}, \quad 
p_3 = \frac{7k^3+6k^2+k}{2}, \quad 
p_4 = \frac{105k^4+154k^3+63k^2+6k}{8},\ldots
$$

Modern computer algebra packages can generally carry out this
type of symbolic summation quickly, so we can use this polynomial
recurrence directly to find $p_n$, and evaluate $p_{n-1}(1)$ to find the
number of nodes on level $n$.  For example, in Maple\tm:
\begin{verbatim}
 p := proc(n) option remember;
      if (n=1) then k else sum( p(n-1), 'k'=k+1..2*k ) fi; end;
 s := n -> eval( p(n-1), k=1 );
\end{verbatim}
Then \verb"p(n)" returns the polynomial $p_n$ in the variable
\verb"k", and \verb"s(n)" evaluates $p_{n-1}$ at $k=1$, giving us the
number of nodes on level $n$ of the tree.

Now that we know that $p_n$ is an $n$th-degree polynomial, we need not
calculate the polynomial explicitly just to find some of its values.
For example, $p_n$ is determined by its values at the $n+1$ points
$k=0,1,\ldots,n$, which we can find (using Recurrence~\ref{rec_d})
once we know $p_{n-1}$ at $k=0,1,\ldots,2n$.  We could then fit an
interpolating polynomial to those points to find other desired values
of $p_n$.  In this case, another trick presents itself; we can use the
linear dependence among $n+2$ equally-spaced values of a polynomial
$p$ of degree $n$:
$$
\forall a,b: \,\,
\sum_{i=0}^{n+1} (-1)^i {n \choose i}\, p(a+bi) = 0.
$$
Combining all of these tricks, we can efficiently calculate the number
of nodes on row $n$ of our tree as follows:
\begin{recur}
\label{rec_good}
\begin{eqnarray*}
d(0,k) &=& 1 \mbox{\quad for all $k$,} \\
d(n,0) &=& 0 \mbox{\quad for all $n>0$,} \\
d(n,k) &=& d(n,k-1) - d(n-1,k) + d(n-1,2k-1) + d(n-1,2k) \\
       && \hspace{2.25in} \mbox{for $k\leq n$, and} \\
d(n,k) &=& \sum_{j=1}^{n+1} (-1)^{(j-1)} {n+1 \choose j} \, d(n,k-j)
       \mbox{\qquad for $n<k\leq 2k+2$.}
\end{eqnarray*}
\end{recur}
As before, $s(n)=d(n-1,1)$ gives the number of nodes on level $n$ of
our tree.

Now that we have a refined computational technique, we would like to
evaluate its efficiency and justify our statement that it can
calculate $s(n)$ in time polynomial in $n$.  The difficulty here is
that we cannot in good conscience simply count the number of
arithmetic operations involved in finding $d(n-1,1)$ using
Recurrence~\ref{rec_good}, because the values we encounter grow
quickly as $n$ increases.  Observe, for example, that $s(n)$ certainly
grows faster than $c^n$ for any fixed constant $c$, because on row $n$
of the tree, every node has at least $n$ (and at most $2^{n-1}$)
children.

In this situation, an analysis of algorithmic complexity must take
into account the magnitude of the numbers involved in arithmetic
operations.  Bach and Shallit~\cite{ant} argue for what they call the
{\em naive bit complexity} measure, in which we can calculate $a+b$
and $ab$ in time $O(\log a + \log b)$ and $O(\log a \log b)$,
respectively.  These time estimates reflect the speed of the naive,
grade-school algorithms for adding and multiplying two numbers with
$\log a$ and $\log b$ digits; the authors argue that these estimates
are both asymptotically realistic and pragmatic for predicting
real-world behavior of computations.

\begin{thm}
\label{thm_poly}
The naive bit complexity of calculating $s(n)$ is $O(n^6)$.
\end{thm}
Note that this is not what one might generally call a polynomial-time
computation, since it is polynomial in $n$, not in $\log n$, the
length of the input.  As already mentioned, $s(n)$ grows so quickly
that it could not even be written down in time polynomial in $\log n$.

To calculate $s(n)$, we will compute all values $d(m,k)$ for $0\leq
m\leq n-1$ and $0\leq k\leq 2m+2$, in lexicographic order, using
Recurrence~\ref{rec_good}.  We first need some bound on the size of
the numbers we encounter.

\begin{lemma}
\label{lem_bound}
%$$
%{\left( {\textstyle\frac32} \right)\!}^{n(n-1)/2}
%\,\leq\,
%\frac{d(n,k)}{k^n}
%\,\leq\,
%2^{n(n-1)/2}
%$$
\quad $d(n,k) \,\leq\, 2^{n(n-1)/2}\, k^n$
\end{lemma}

\pf
The bound is based on Recurrence~\ref{rec_d}, $d(n,k) =
\sum_{j=k+1}^{2k} d(n-1,j)$.  For a fixed $n\geq1$, we know $d(n,k)$
is an increasing function on positive integers $k$.  Therefore the largest
of the $k$ terms in the sum is $d(n-1,2k)$, and we have $d(n,k)\leq
k\,d(n-1,2k)$.  Repeating, we get:
\begin{eqnarray*}
d(n,k) &\leq& k\,d(n-1,2k) \,\,\leq\,\, (k)(2k)\,d(n-2,4k) \,\,\leq \cdots \\
       &\leq& (k)(2k)\cdots(2^{n-1}k)\,d(0,2^nk) \\
       &    & =\,\, 2^{n(n-1)/2} \, k^n 
\end{eqnarray*}
In the next section we will discuss the asymptotic growth of $s(n)$,
which we have just bounded by $2^{n-1 \choose 2}$.
%which we have just bounded by $2^{(n-1)(n-2)/2}.$
\qed

\noindent{\bf Proof\,} of Theorem~\ref{thm_poly}:
We will perform the calculation using Recurrence~\ref{rec_good},
storing partial results so we never calculate any value of $d$ twice.

Fix some $n$, and suppose we know $d(n-1,k)$ for all $0\leq k\leq 2n$.
Calculating $d(n,k)$ for all $0\leq k\leq 2n+2$ involves two phases.
\begin{enumerate}
\item
For each $k$ with $0\leq k\leq n$, we compute the sum of four numbers.
By Lemma~\ref{lem_bound}, the summands are of length $O(n^2 + n\log k)$, 
which is $O(n^2)$ since $k$ is small.  Thus each value of $k$ takes
time $O(n^2)$, and the whole phase takes time $O(n^3)$.

\item
For $n+1\leq k\leq 2n+2$, we compute $\sum_{j=1}^{n+1} (-1)^{(j-1)}
{n+1 \choose j} d(n,k-j)$.  As above, $d(n,k-j)$ has length $O(n^2)$.
The binomial coefficient ${n+1 \choose j}$ has length $O(n)$, since
numbers in the $n$th row of Pascal's Triangle are bounded by $2^n$;
the cost of computing it is negligible because we are using all of the
top $n$ rows of Pascal's Triangle, which we can compute by additive
recurrence.

Then we form each product in time $O(n^3)$ and take their
sum in time $O(n^4)$ for each $k$.  Thus the whole phase takes time
$O(n^5)$.
\end{enumerate}
Thus passing from $n-1$ to $n$ takes time $O(n^5)$, and we can
calculate $s(n)$ from scratch in time $O(n^6)$.
\qed

In practice, Recurrence~\ref{rec_good} is easy to implement and
seems to perform much better than the above analysis suggests for
values of $n$ we are interested in.  In Maple\tm:
\begin{verbatim}
d := proc(n,k) 
  local j;
  option remember;
  if n=0 then 1
    elif k=0 then 0
    elif k<=n then d(n,k-1)-d(n-1,k)+d(n-1,2*k-1)+d(n-1,2*k)
    else add( (-1)^(j-1) * binomial(n+1,j) * d(n,k-j), j=1..n+1 )
  fi
end;
s := n -> d(n-1,1);
\end{verbatim}
This code seems empirically to calculate $s(n)$ after having done the
work for $s(n-1)$ in time $O(n^2)$ even when $n$ is around 190, when
$s(n)$ has over 5000 decimal digits.  This is the behavior one would
expect if multiplication had a constant unit cost and addition were
free.

On a modest desktop Pentium II, this computes up to $s(30)$ in under a
second, $s(85)$ in under a minute, and s(190) in about an hour; a
little extra work to avoid computing the binomial coefficients
multiple times speeds it up even more.  We record $s(n)$ for $1\leq
n\leq 22$ here.  The sequence also appears as entry A008934 in
Sloane's On-Line Encyclopedia of Integer Sequences~\cite{EIS}.

\begin{quote} \flushleft \sloppy
1, 1, 2, 7, 41, 397, 6377, 171886, 7892642, 627340987, 87635138366,
   21808110976027, 9780286524758582, 7981750158298108606,
   11950197013167283686587, 33046443615914736611839942, 
   169758733825407174485685959261,
   1627880269212042994531083889564192, 
   29264239787495935863325877024506142042,
   989901541366810465070950556260422637919176,
   63214893835996484808167529681187283166038800097,
   7643667309922877343580868981767361594845888953165967,
   \ldots
\end{quote}

\subsection{Asymptotic Behavior}

Now we turn to the asymptotic growth of $s(n)$.  We know from
Lemma~\ref{lem_bound} that ${n-1 \choose 2}$ is an upper bound for
$\lg s(n)$, where $\lg$ denotes $\log_2$.  We will first prove that a
lower bound for $s(n)$ is $\alpha\,2^{n\choose2} / (n-1)!$ for a
certain constant $\alpha$, and then we will show that $\lg s(n)$ is
asymptotic to ${n\choose2}-\lg n!+O(\log n)^2$.

We are grateful to Donald Knuth for suggesting the method used here.
The technique rests on the following observation:

\begin{thm}[Knuth]
\label{thm_knuth}
Let $T$ be a rooted tree in which we want to know the number of
vertices on the $n$th level.  Select $n-1$ vertices
$v_1,\ldots,v_{n-1}$ in $T$ by choosing $v_1$ to be the root and
picking $v_{i+1}$ uniformly at random from among the children of
$v_i$, of which there are $\deg(v_i)$.  Then the expected value
$$
E(\deg(v_1)\deg(v_2)\cdots\deg(v_{n-1}))
$$
is exactly the number of vertices on the $n$th level of $T$.
\end{thm}

\pf
See \cite{knuth} for a discussion of this technique in greater
generality.  In this case, a proof by induction is straightforward: if
this technique works for counting the $n$th level of each of $k$ trees
$T_1,\ldots,T_k$, then it clearly also works for counting level $n+1$
of the tree $T$ whose root has $k$ children, the roots of
$T_1,\ldots,T_k$.
\qed

We will apply Theorem~\ref{thm_knuth} to the tree of tournament
sequences, in which, conveniently, each vertex is already labelled
with its degree.  This means we can calculate $s(n)$ by finding the
expected value of the product $t_1t_2\cdots t_{n-1}$, where
$(t_1,t_2,\ldots,t_{n-1})$ is a tournament sequence selected at random
by setting $t_1=1$ and picking $t_{i+1}$ uniformly at random from
among $t_i+1,t_i+2,\ldots,2t_i$.  For the remainder of this section,
whenever we talk about a distribution for $t_i$, it is implicitly with
respect to this way of picking a tournament sequence.

\begin{lemma}
\label{lem_lb}
Let $s(n)$ be the number of tournament sequences of length $n$.  Then
$$
s(n) \geq \alpha\,\frac{2^{n\choose 2}}{(n-1)!}
$$
where $\alpha=(1-\frac12)(1-\frac14)(1-\frac18)\cdots\approx.28878837\ldots$
\end{lemma}

\pf
Observe that there is a natural continuous analogue to the expected
value $E(t_1t_2\cdots t_{n-1})=s(n)$.  Consider instead the expected
value $E(r_1r_2\cdots r_{n-1})$, where $r_1=1$ and $r_{i+1}$ is a real
number chosen uniformly at random from the interval $(r_i,2r_i]$.
Equivalently, we are taking random variables $u_i$ $(1\leq i\leq n-2)$
each with a uniform distribution over $(1,2]$, and setting
$r_{i+1}=u_ir_i$.

The expected value in the continuous case will give an underestimate
of the expected value for the discrete version, in which the $u_i$ are
distributed the same way but we set $t_{i+1}=\lceil u_it_i \rceil$.
The independence of the various $u_i$ make the expected value easy to
calculate:
\begin{eqnarray*}
  E(r_1r_2r_3\cdots r_{n-1})
  &=& E((1)(u_1)(u_1u_2)\cdots(u_1u_2\cdots u_{n-2})) \\
  &=& E(u_1^{n-2} \, u_2^{n-3} \cdots u_{n-2}^1) \\
  &=& \frac{2^{n-1}-1}{n-1}\,\frac{2^{n-2}-1}{n-2}\cdots\frac{2^2-1}{2} \\
&\geq& \alpha\,\frac{2^{n\choose 2}}{(n-1)!}
\end{eqnarray*}
\qed

Now we show that this lower bound is quite good, by finding the rate
of growth of the error.  We use $\lg$ to denote $\log_2$.

\begin{thm}
\label{thm_growth}
\quad $\lg s(n) = {n \choose 2} - \lg n! + O(\log n)^2$
\end{thm}

\pf
In the proof of Lemma~\ref{lem_lb}, we calculated $E(u_k^{n-k-1})$,
where $u_k=r_{k+1}/r_k$ was uniformly distributed over $(1,2]$.  In
the original problem, $t_k$ and $t_{k+1}$ are integers, so the ratio
$t_{k+1}/t_k$ instead takes on a discrete set of values, each with
equal probability:
$$
E\left( \left(\frac{t_{k+1}}{t_k}\right)^{\!\!p} \,\right) = 
 \frac{ \left(\frac{t+1}{t}\right)^p +
        \left(\frac{t+2}{t}\right)^p + \cdots +
        \left(\frac{t+t}{t}\right)^p }{t}
$$
where $t=t_k$ throughout.  We can rewrite this and take advantage of
the ability to sum consecutive $p$th powers:
\begin{eqnarray*}
E\left( \left(\frac{t_{k+1}}{t_k}\right)^{\!\!p} \,\right)
&=&
  \frac{ (t+1)^p + (t+2)^p +\cdots+ (t+t)^p }{t^{p+1}} \\
&\leq& 
  \frac{\frac{(2t)^{p+1}}{p+1} + \frac{(2t)^p}{2} + O(2t)^{p-1}}{t^{p+1}} \\
&&=\,
  \frac{2^{p+1}}{p+1} + \frac{2^p}{2t} + O\!\left(\frac{2^p}{t^2}\right)
\end{eqnarray*}
The $2^{p+1}/p+1$ term is exactly the lower bound we used in
Lemma~\ref{lem_lb}, and we now have an idea of how much error this
introduced:
\begin{eqnarray*}
E\left( \left(\frac{t_{k+1}}{t_k}\right)^{\!\!n-k-1} \right)
&=&
\frac{2^{n-k}}{n-k} + \frac{2^{n-k}}{4t_k} 
                    + O\!\left(\frac{2^{n-k}}{t_k^2}\right)
\\ &=& 
E(u_k^{n-k-1}) \left( 1 + \frac{n-k}{4t_k} 
                        + O\!\left(\frac{n-k}{t_k^2}\right) \right)
\end{eqnarray*}
The expected value now depends on $t_k$, and to bound the error, we
need some idea of how large we expect $t_k$ to be.

Consider again the continuous analogue used in the proof of
Lemma~\ref{lem_lb}.  We expect $s_k$ to grow exponentially, as $c^k$
for some $c$, so look at the distribution of $s_k^{1/k}$.  Taking
logs, we see that $\log s_k^{1/k}$ is distributed as the average of
$k$ copies of $\log x$ on $(1,2]$, each with mean $2\log2-1$.  Thus as
$k$ increases without bound, the distribution of $s_k^{1/k}$ converges
towards a point distribution at $4/e$.  As we already noted, the $t_k$
certainly grow no slower than the $s_k$.  We conclude that for any
constant $c<4/e$ and probability $p<1$, there is a sufficiently
large $k_0$ such that $\Pr(t_k>c^k)>p$ for all $k\geq k_0$.

The error we want to bound is the product of the $n$ error terms $e_k
= 1+\frac{n-k}{4t_k}+O(\frac{n-k}{t_k^2})$ for $k=1,\ldots,n$.  To
simplify bookkeeping, we will instead think about $\lg s(n)$ and bound
the sum of the logs of the error terms.  We separate the work into two
cases, depending on whether $k$ is less or greater than $2\log n$.

When $k<2\log n$, the denominator $t_k$ may be small, and the error
terms $\log e_k$ may be as much as $O(\log n)$.  Adding up all $2\log
n$ of these terms gives a total which is $O(\log n)^2$, the error term
in the statement of our theorem.

Finally, when $k=2\log n$, we know that as long as $n$ is sufficiently
large, with high probability $t_k>c^{2\log n}$, which is on the order
of $n^2$.  Thus with high probability, $e_k$ is $1+O(1/n)$, and $\log
e_k=O(1/n)$.  Since the error terms monotonically decrease, the sum of
all the terms with $k\geq 2\log n$ is bounded by $ne_{2\log n}=O(1)$.
So taking $n$ sufficiently large ensures that the error in the lower
bound is concentrated almost entirely in the first $2\log n$ error
terms, and we are done.
\qed

Computational evidence based on the actual values of $s(n)$ for $n$ up
to 190 indicates that the constant needed to make $\lg s(n) <
{n\choose 2} - \lg n! + c (\log n)^2$ reaches a peak of $c\approx
1.18304060\ldots$ at $n=32$ and decreases slowly thereafter.

\end{document}